\documentclass[12pt]{article}
\usepackage{color}
\begin{document}

 \begin{center}
   {\large  A note on log-convexity of $q$-Catalan numbers}\\*[5pt]
    L.M. Butler and W.P. Flanigan\footnotemark
 \end{center}

 \footnotetext{Research initiated in Pat Flanigan's senior thesis, under Lynne Butler's direction at Haverford College, Haverford PA 19041 USA.}

\bigskip

The sequence of Catalan numbers
$(1, 1, 2, 5, 14, \ldots)$
is defined by $C_0=1$ and the recursion 
\[ C_{n+1}=\sum_{k=0}^n C_kC_{n-k},\]
from which it is evident that the sequence is nondecreasing. From the well-known formula 
$C_n={1\over n+1}\left( {2n\atop n}\right)$
it is also easily seen that
$C_{k-1}C_{k+1}\geq {C_k}^2$ for $k\geq 1$.
Equivalently, the sequence satisfies the log-convexity result that 
\[ C_{k-1}C_{\ell+1}-C_kC_\ell \mbox{ is nonnegative for $\ell\geq k\geq 1$}.\]

Many polynomials in $q$  have been called $q$-Catalan numbers. See [3]. These polynomials have nonnegative coefficients and evaluate at $q=1$ to  Catalan numbers. One such sequence,  studied by  Carlitz and Riordan in [2],
$ (1, 1, 1+q, 1+q+2q^2+q^3, 1+q+2q^2+3q^3+3q^4+3q^5+q^6, \ldots)$
is defined by  
$C_0(q)=1$ and the recursion
\[ C_{n+1}(q)=\sum_{k=0}^n q^{(k+1)(n-k)}C_k(q)C_{n-k}(q).\]
The recursion shows that the sequence is nondecreasing: $C_{n+1}(q)-C_n(q)$ has nonnegative coefficients for $n\geq 0$. No formula is known for $C_n(q)$ and it is not true that $C_2(q)C_4(q)-{C_3(q)}^2$ has nonnegative coefficients. But it is true that  $C_{k-1}(q)C_{k+1}(q)-q{C_k(q)}^2$ has nonnegative coefficients for $k\geq 1$. Our main result that 
\[ C_{k-1}(q)C_{\ell+1}(q)-q^{\ell-k+1}C_k(q)C_\ell(q) \mbox{ has nonnegative coefficients for $\ell\geq k\geq 1$}\]
is more general than the special case for $k=\ell$. 

\medskip

Other authors [4] have defined a sequence of polynomials $\{P_n(q)\}_{n\geq 0}$ to be $q$-log-convex if $P_{k-1}(q)P_{k+1}(q)-{P_k(q)}^2$ has nonnegative coefficients for $k\geq 1$. We do not advocate this definition, since it does not imply that $P_{k-1}(q)P_{\ell+1}(q)-P_k(q)P_\ell(q)$ has nonnegative coefficients for $\ell\geq k\geq 1$. Consider $1+q+q^2+q^5, 1+q^2+q^3,1+q+q^3,1+2q+q^6$. Log-convexity for sequences of nonzero polynomials with nonnegative coefficients is more subtle than log-convexity for sequences of positive numbers.

To prove our log-convexity result, we use the combinatorial interpretation 
\[ C_k(q)=\sum_\pi q^{{\rm inv}\, \pi},\] 
where the sum is over   permutations with  $k$ 1s and $k$ 2s such that every initial segment has no more 2s than 1s. For example, $C_3(q)=1+q+2q^2+q^3$, since the permutations 111222, 112122, 121122, 112212  and 121212  have inversion numbers 0, 1, 2, 2 and 3, respectively. Such a permutation may be visualized as a path in a $k\times k$ square from the lower left corner to the upper right corner. If the permutation is read left to right, where 1 denotes a horizontal step and 2 denotes a vertical step, then the path never rises above the diagonal. The inversion number of the permutation is the area below the path. Such permutations, called lattice permutations, were studied by MacMahon [5].

\medskip

\noindent {\bf Theorem}: The $q$-Catalan numbers $C_k(q)=\sum_\pi q^{{\rm inv}\, \pi}$, where the sum is over lattice permutations with $k$ 1s and $k$ 2s,  satisfy: \[C_{k-1}(q)C_{\ell+1}(q)-q^{\ell-k+1}C_k(q)C_\ell(q) \;\mbox{has nonnegative coefficients for $k\leq \ell$.}\]

\medskip

The proof of the above result uses an injection introduced by Butler [1] to prove a similar log-concavity result for $q$-binomial coefficients: 

\[ \left[ {n\atop k}\right]_q \left[ {n\atop \ell}\right]_q -q^{\ell-k+1}\left[ {n\atop k-1}\right]_q \left[ {n\atop \ell+1}\right]_q  \;\mbox{has nonnegative coefficients for $k\leq \ell$.}\]

\medskip

\noindent {\bf Proof of the Theorem}: Given $k\leq \ell$, we define an injection
\[ \varphi: {\cal P}_k\times {\cal P}_\ell\rightarrow {\cal P}_{k-1}\times {\cal P}_{\ell+1}\]
where ${\cal P}_n$ is the set of lattice permutations  with  $n$ 1s and $n$ 2s. Then we show that if $\varphi(\pi,\sigma)=(\nu,\omega)$, we have ${\rm inv}\,\pi+{\rm inv}\,\sigma+(\ell-k+1)={\rm inv}\,\nu+{\rm inv}\,\omega$. 

\medskip

Let $(\pi,\sigma)\in {\cal P}_k\times {\cal P}_\ell$ be given. For $0\leq i\leq 2k-2$, consider
\begin{eqnarray*}
\nu^{(i)}&=&\sigma_1\cdots\sigma_i\pi_{i+3}\cdots\pi_{2k}\\
\omega^{(i)}&=&\pi_1\cdots\pi_i\pi_{i+1}\pi_{i+2}\sigma_{i+1}\cdots\sigma_{2\ell}\end{eqnarray*}
Define $\varphi(\pi,\sigma)=(\nu^{(i)},\omega^{(i)})=(\sigma_L\pi_R,\pi_L\sigma_R)$, where $i$ is smallest such that the number of 2s in $\pi_L$ exceeds the number of 2s in $\sigma_L$ by exactly 1. (Note that $\pi_L=\pi_1\pi_2$ and $\sigma_L=\emptyset$ if $i=0$, but $\pi_L=\pi$ and $\sigma_L=\sigma_1\cdots \sigma_{2k-2}$ if $i=2k-2$. Initially $\pi_L$ has at most one more 2 than $\sigma_L$, and finally $\pi_L$ has at least one more 2 than $\sigma_L$.)

\medskip

To show that $(\ell-k+1)+{\rm inv}\,\pi+{\rm inv}\,\sigma={\rm inv}\,\sigma_L\pi_R+{\rm inv}\,\pi_L\sigma_R$, we use the fact that the number of 2s in $\pi_L$, denoted $m_2\pi_L$, equals $m_2\sigma_L+1$; hence the number of 1s in $\pi_L$, denoted $m_1\pi_L$, equals $m_1\sigma_L+1$. Since an inversion 21 in a permutation may occur in the left portion, occur in the right portion, or straddle the left and right portions, we have:
\begin{eqnarray*}
{\rm inv}\,\pi&=&{\rm inv}\,\pi_L+{\rm inv}\,\pi_R+(m_2\pi_L)(m_1\pi_R)\\
{\rm inv}\,\sigma&=&{\rm inv}\,\sigma_L+{\rm inv}\,\sigma_R+(m_2\sigma_L)(m_1\sigma_R)\\
{\rm inv}\,\pi_L\sigma_R&=&{\rm inv}\,\pi_L+{\rm inv}\,\sigma_R+(m_2\pi_L)(m_1\sigma_R)\\
{\rm inv}\,\sigma_L\pi_R&=&{\rm inv}\,\sigma_L+{\rm inv}\,\pi_R+(m_2\sigma_L)(m_1\pi_R)\end{eqnarray*}
\begin{eqnarray*}
{\rm inv}\,\pi_L\sigma_R+{\rm inv}\,\sigma_L\pi_R-{\rm inv}\,\pi-{\rm inv}\,\sigma&=&(m_2\pi_L-m_2\sigma_L)(m_1\sigma_R-m_1\pi_R)\\
&=&m_1\sigma_R-m_1\pi_R\\
&=&(\ell-m_1\sigma_L)-(k-m_1\pi_L)\\
&=&\ell-k+1\end{eqnarray*}

\bigskip

\noindent {\bf Corollary:} For $1\leq r \leq k$ and $\ell > k-r$,
 \[C_{k-r}(q)C_{\ell+r}(q)-q^{r(\ell-k+r)}C_k(q)C_\ell(q) \;\mbox{has nonnegative coefficients.}\]
 
\noindent Note that no higher power of $q$ may be used. To see this, note that $C_n(q)$ is monic of degree $\left( {n\atop 2}\right)$.  So ${\rm deg}\,C_{k-r}(q)C_{\ell+r}(q)={\rm deg}\,q^{r(\ell-k+r)}C_k(q)C_\ell(q).$

\medskip

This corollary may be proved by induction on $r$ or may be proved using an injection like the one above  for $r=1$. To visualize this injection $\varphi$, picture $\pi$ as a lattice path $L_\pi$ from the lower left corner to the upper right corner of a $k\times k$ square and picture $\sigma$ as a lattice path $L_\sigma$ from the lower left corner to the upper right corner of a $\ell\times \ell$ square. Arrange these inside a $(\ell+r)\times (\ell +r)$ square with the lattice path $L_\pi$ beginning in the lower left corner and the lattice path $L_\sigma$ ending in the upper right corner. Find the point where $L_\pi$ first meets $L_\sigma$. If $\varphi(\pi,\sigma)=(\nu,\omega)$, then the lattice path $L_\omega$ follows $L_\pi$ until this point then follows $L_\sigma$, and the lattice path $L_\nu$ follows $L_\sigma$ until this point then follows $L_\pi$. The exponent $r(\ell-k+r)$ is the area of the rectangle in the lower right of the $(\ell+r)\times (\ell+r)$ square that lies to the right of the $k\times k$ square in the lower left and lies below the $\ell\times \ell $ square in the upper right. 

\newpage

\noindent {\bf Example:} Terms in  
$ q^{2(3)}\textcolor{red}{C_6(q)}\textcolor{blue}{C_7(q) }$ and $ \textcolor{black}{C_{4}(q)}C_{9}(q)$, are compared below:\newline
\hspace*{1.1in}$ \textcolor{red}{112112221122} \quad\quad\quad\;\;\;{}_{{}_{{}_{\displaystyle \mapsto}}}\quad\quad \textcolor{blue}{ 121}\textcolor{red}{21122}$\newline
\hspace*{1.1in}$\phantom{1121}\textcolor{blue}{    12111212212212}\qquad \textcolor{red}{1121122}\textcolor{blue}{ 11212212212}$
\[
\setlength{\unitlength}{11pt}
\begin{picture}(0,9)(6.05,0)
\multiput(2.25,0)(1,0){4}{\textcolor{red}{\line(1,0){.5}}}
\multiput(4.25,1)(1,0){2}{\textcolor{red}{\line(1,0){.5}}}
\multiput(4.25,2)(1,0){2}{\textcolor{red}{\line(1,0){.5}}}
\multiput(3,.25)(0,1){1}{\textcolor{red}{\line(0,1){.5}}}
\multiput(4,.25)(0,1){1}{\textcolor{red}{\line(0,1){.5}}}
\multiput(5,.25)(0,1){4}{\textcolor{red}{\line(0,1){.5}}}
\multiput(6,.25)(0,1){3}{\textcolor{red}{\line(0,1){.5}}}
\multiput(3.25,2)(1,0){1}{\textcolor{blue}{\line(1,0){.5}}}
\multiput(6.25,2)(1,0){3}{\textcolor{blue}{\line(1,0){.5}}}
\multiput(6.25,3)(1,0){3}{\textcolor{blue}{\line(1,0){.5}}}
\multiput(7.25,4)(1,0){2}{\textcolor{blue}{\line(1,0){.5}}}
\multiput(7.25,5)(1,0){2}{\textcolor{blue}{\line(1,0){.5}}}
\multiput(8.25,6)(1,0){1}{\textcolor{blue}{\line(1,0){.5}}}
\multiput(8.25,7)(1,0){1}{\textcolor{blue}{\line(1,0){.5}}}
\multiput(7,2.25)(0,1){2}{\textcolor{blue}{\line(0,1){.5}}}
\multiput(8,2.25)(0,1){4}{\textcolor{blue}{\line(0,1){.5}}}
\multiput(9,2.25)(0,1){6}{\textcolor{blue}{\line(0,1){.5}}}
\thicklines
\put(0,0){\textcolor{red}{\line(1,0){2}}}
\put(2,0){\textcolor{red}{\line(0,1){1}}}
\put(2,1){\textcolor{red}{\line(1,0){2}}}
\put(4,1){\textcolor{red}{\line(0,1){3}}}
\put(4,4){\textcolor{red}{\line(1,0){2}}}
\put(6,4){\textcolor{red}{\line(0,1){2}}}
\put(2,2){\textcolor{blue}{\line(1,0){1}}}
\put(3,2){\textcolor{blue}{\line(0,1){1}}}
\put(3,3){\textcolor{blue}{\line(1,0){3}}}
\put(6,3){\textcolor{blue}{\line(0,1){1}}}
\put(6,4){\textcolor{blue}{\line(1,0){1}}}
\put(7,4){\textcolor{blue}{\line(0,1){2}}}
\put(7,6){\textcolor{blue}{\line(1,0){1}}}
\put(8,6){\textcolor{blue}{\line(0,1){2}}}
\put(8,8){\textcolor{blue}{\line(1,0){1}}}
\put(9,8){\textcolor{blue}{\line(0,1){1}}}
\thinlines
\put(12,0){
\begin{picture}(9,9)
\multiput(2.25,0)(1,0){4}{\textcolor{red}{\line(1,0){.5}}}
\multiput(4.25,1)(1,0){2}{\textcolor{red}{\line(1,0){.5}}}
\multiput(4.25,2)(1,0){2}{\textcolor{red}{\line(1,0){.5}}}
\multiput(3,.25)(0,1){1}{\textcolor{red}{\line(0,1){.5}}}
\multiput(4,.25)(0,1){1}{\textcolor{red}{\line(0,1){.5}}}
\multiput(5,.25)(0,1){4}{\textcolor{red}{\line(0,1){.5}}}
\multiput(6,.25)(0,1){3}{\textcolor{red}{\line(0,1){.5}}}
\multiput(3.25,2)(1,0){1}{\textcolor{blue}{\line(1,0){.5}}}
\multiput(6.25,2)(1,0){3}{\textcolor{blue}{\line(1,0){.5}}}
\multiput(6.25,3)(1,0){3}{\textcolor{blue}{\line(1,0){.5}}}
\multiput(7.25,4)(1,0){2}{\textcolor{blue}{\line(1,0){.5}}}
\multiput(7.25,5)(1,0){2}{\textcolor{blue}{\line(1,0){.5}}}
\multiput(8.25,6)(1,0){1}{\textcolor{blue}{\line(1,0){.5}}}
\multiput(8.25,7)(1,0){1}{\textcolor{blue}{\line(1,0){.5}}}
\multiput(7,2.25)(0,1){2}{\textcolor{blue}{\line(0,1){.5}}}
\multiput(8,2.25)(0,1){4}{\textcolor{blue}{\line(0,1){.5}}}
\multiput(9,2.25)(0,1){6}{\textcolor{blue}{\line(0,1){.5}}}
\multiput(6.25,0)(1,0){3}{\line(1,0){.5}}
\multiput(6.25,1)(1,0){3}{\line(1,0){.5}}
\multiput(7,.25)(0,1){2}{\line(0,1){.5}}
\multiput(8,.25)(0,1){2}{\line(0,1){.5}}
\multiput(9,.25)(0,1){2}{\line(0,1){.5}}
\thicklines
\put(0,0){\textcolor{red}{\line(1,0){2}}}
\put(2,0){\textcolor{red}{\line(0,1){1}}}
\put(2,1){\textcolor{red}{\line(1,0){2}}}
\put(4,1){\textcolor{red}{\line(0,1){3}}}
\put(4,4){\textcolor{red}{\line(1,0){2}}}
\put(6,4){\textcolor{red}{\line(0,1){2}}}
\put(2,2){\textcolor{blue}{\line(1,0){1}}}
\put(3,2){\textcolor{blue}{\line(0,1){1}}}
\put(3,3){\textcolor{blue}{\line(1,0){3}}}
\put(6,3){\textcolor{blue}{\line(0,1){1}}}
\put(6,4){\textcolor{blue}{\line(1,0){1}}}
\put(7,4){\textcolor{blue}{\line(0,1){2}}}
\put(7,6){\textcolor{blue}{\line(1,0){1}}}
\put(8,6){\textcolor{blue}{\line(0,1){2}}}
\put(8,8){\textcolor{blue}{\line(1,0){1}}}
\put(9,8){\textcolor{blue}{\line(0,1){1}}}
\end{picture}}
\end{picture}
\begin{picture}(9,9)(6,0)
\multiput(2.25,0)(1,0){4}{\textcolor{red}{\line(1,0){.5}}}
\multiput(4.25,1)(1,0){2}{\textcolor{red}{\line(1,0){.5}}}
\multiput(4.25,2)(1,0){2}{\textcolor{red}{\line(1,0){.5}}}
\multiput(3,.25)(0,1){1}{\textcolor{red}{\line(0,1){.5}}}
\multiput(4,.25)(0,1){1}{\textcolor{red}{\line(0,1){.5}}}
\multiput(5,.25)(0,1){4}{\textcolor{red}{\line(0,1){.5}}}
\multiput(6,.25)(0,1){3}{\textcolor{red}{\line(0,1){.5}}}
\multiput(3.25,2)(1,0){1}{\textcolor{blue}{\line(1,0){.5}}}
\multiput(6.25,2)(1,0){3}{\textcolor{blue}{\line(1,0){.5}}}
\multiput(6.25,3)(1,0){3}{\textcolor{blue}{\line(1,0){.5}}}
\multiput(7.25,4)(1,0){2}{\textcolor{blue}{\line(1,0){.5}}}
\multiput(7.25,5)(1,0){2}{\textcolor{blue}{\line(1,0){.5}}}
\multiput(8.25,6)(1,0){1}{\textcolor{blue}{\line(1,0){.5}}}
\multiput(8.25,7)(1,0){1}{\textcolor{blue}{\line(1,0){.5}}}
\multiput(7,2.25)(0,1){2}{\textcolor{blue}{\line(0,1){.5}}}
\multiput(8,2.25)(0,1){4}{\textcolor{blue}{\line(0,1){.5}}}
\multiput(9,2.25)(0,1){6}{\textcolor{blue}{\line(0,1){.5}}}
\put(0,0){\line(0,1){6}}
\put(0,6){\line(1,0){6}}
\put(2,2){\line(0,1){7}}
\put(2,9){\line(1,0){7}}
\thicklines
\put(0,0){\textcolor{red}{\line(1,0){2}}}
\put(2,0){\textcolor{red}{\line(0,1){1}}}
\put(2,1){\textcolor{red}{\line(1,0){2}}}
\put(4,1){\textcolor{red}{\line(0,1){3}}}
\put(4,4){\textcolor{red}{\line(1,0){2}}}
\put(6,4){\textcolor{red}{\line(0,1){2}}}
\put(2,2){\textcolor{blue}{\line(1,0){1}}}
\put(3,2){\textcolor{blue}{\line(0,1){1}}}
\put(3,3){\textcolor{blue}{\line(1,0){3}}}
\put(6,3){\textcolor{blue}{\line(0,1){1}}}
\put(6,4){\textcolor{blue}{\line(1,0){1}}}
\put(7,4){\textcolor{blue}{\line(0,1){2}}}
\put(7,6){\textcolor{blue}{\line(1,0){1}}}
\put(8,6){\textcolor{blue}{\line(0,1){2}}}
\put(8,8){\textcolor{blue}{\line(1,0){1}}}
\put(9,8){\textcolor{blue}{\line(0,1){1}}}
\put(0,0){\circle{.3}}
\put(9,9){\circle*{.3}}
\put(2,2){\circle*{.3}}
\put(6,6){\circle{.3}}
\put(10.5,4){$\mapsto$}
\thinlines
\put(12,0){
\begin{picture}(9,9)
\multiput(2.25,0)(1,0){4}{\textcolor{red}{\line(1,0){.5}}}
\multiput(4.25,1)(1,0){2}{\textcolor{red}{\line(1,0){.5}}}
\multiput(4.25,2)(1,0){2}{\textcolor{red}{\line(1,0){.5}}}
\multiput(3,.25)(0,1){1}{\textcolor{red}{\line(0,1){.5}}}
\multiput(4,.25)(0,1){1}{\textcolor{red}{\line(0,1){.5}}}
\multiput(5,.25)(0,1){4}{\textcolor{red}{\line(0,1){.5}}}
\multiput(6,.25)(0,1){3}{\textcolor{red}{\line(0,1){.5}}}
\multiput(3.25,2)(1,0){1}{\textcolor{blue}{\line(1,0){.5}}}
\multiput(6.25,2)(1,0){3}{\textcolor{blue}{\line(1,0){.5}}}
\multiput(6.25,3)(1,0){3}{\textcolor{blue}{\line(1,0){.5}}}
\multiput(7.25,4)(1,0){2}{\textcolor{blue}{\line(1,0){.5}}}
\multiput(7.25,5)(1,0){2}{\textcolor{blue}{\line(1,0){.5}}}
\multiput(8.25,6)(1,0){1}{\textcolor{blue}{\line(1,0){.5}}}
\multiput(8.25,7)(1,0){1}{\textcolor{blue}{\line(1,0){.5}}}
\multiput(7,2.25)(0,1){2}{\textcolor{blue}{\line(0,1){.5}}}
\multiput(8,2.25)(0,1){4}{\textcolor{blue}{\line(0,1){.5}}}
\multiput(9,2.25)(0,1){6}{\textcolor{blue}{\line(0,1){.5}}}
\multiput(6.25,0)(1,0){3}{\line(1,0){.5}}
\multiput(6.25,1)(1,0){3}{\line(1,0){.5}}
\multiput(7,.25)(0,1){2}{\line(0,1){.5}}
\multiput(8,.25)(0,1){2}{\line(0,1){.5}}
\multiput(9,.25)(0,1){2}{\line(0,1){.5}}
\put(0,9){\line(1,0){9}}
\put(0,0){\line(0,1){9}}
\put(2,2){\line(0,1){4}}
\put(2,6){\line(1,0){4}}
\thicklines
\put(0,0){\textcolor{red}{\line(1,0){2}}}
\put(2,0){\textcolor{red}{\line(0,1){1}}}
\put(2,1){\textcolor{red}{\line(1,0){2}}}
\put(4,1){\textcolor{red}{\line(0,1){3}}}
\put(4,4){\textcolor{red}{\line(1,0){2}}}
\put(6,4){\textcolor{red}{\line(0,1){2}}}
\put(2,2){\textcolor{blue}{\line(1,0){1}}}
\put(3,2){\textcolor{blue}{\line(0,1){1}}}
\put(3,3){\textcolor{blue}{\line(1,0){3}}}
\put(6,3){\textcolor{blue}{\line(0,1){1}}}
\put(6,4){\textcolor{blue}{\line(1,0){1}}}
\put(7,4){\textcolor{blue}{\line(0,1){2}}}
\put(7,6){\textcolor{blue}{\line(1,0){1}}}
\put(8,6){\textcolor{blue}{\line(0,1){2}}}
\put(8,8){\textcolor{blue}{\line(1,0){1}}}
\put(9,8){\textcolor{blue}{\line(0,1){1}}}
\put(0,0){\circle*{.3}}
\put(9,9){\circle*{.3}}
\put(2,2){\circle{.3}}
\put(6,6){\circle{.3}}
\end{picture}
}
\end{picture}
\]
The example shows how the injection from $\textcolor{red}{ {\cal P}_6}\times
 \textcolor{blue}{ {\cal P}_7}$ to $\textcolor{black}{ {\cal P}_4}\times {\cal P}_9$ maps 
 \begin{eqnarray*}
 (\textcolor{red}{\pi},\textcolor{blue}{\sigma})&=&(\textcolor{red}{112112221122},\textcolor{blue}{    12111212212212})\mbox{ to }\\
 (\textcolor{blue}{\sigma_L}\textcolor{red}{\pi_R},\textcolor{red} {\pi_L}\textcolor{blue}{\sigma_R})&=&(\textcolor{blue}{ 121}\textcolor{red}{21122}, \textcolor{red}{1121122}\textcolor{blue}{ 11212212212}).\end{eqnarray*}
 The permutation $\textcolor{red}{\pi}$ is the red lattice path in the $6\times 6$ square, and  $\textcolor{blue}{\sigma}$ is the blue lattice path in the $7\times 7$ square. The $6\times 6$ and $7\times 7$ square overlap inside the $9\times 9$ square at left. At right, a $4\times 4$ square is inside the $9\times 9$ square. In the $4\times 4$ square is the lattice path $\textcolor{blue}{\sigma_L}\textcolor{red}{\pi_R}$, and in the $9\times 9$ square is the lattice path $\textcolor{red} {\pi_L}\textcolor{blue}{\sigma_R}$. The term in $ q^{2(3)}\textcolor{red}{C_6(q)}\textcolor{blue}{C_7(q) }$ corresponding to $ (\textcolor{red}{\pi},\textcolor{blue}{\sigma})$ is $q^{2(3)+\textcolor{red}{10}+\textcolor{blue}{15}}$, because the area under $\textcolor{red}{\pi}$ is $\textcolor{red}{10}$ and the area under  $\textcolor{blue}{\sigma}$ is $\textcolor{blue}{15}$. The term in $ \textcolor{black}{C_{4}(q)}C_{9}(q)$ corresponding to $ (\textcolor{blue}{\sigma_L}\textcolor{red}{\pi_R},\textcolor{red} {\pi_L}\textcolor{blue}{\sigma_R})$ is $q^{5+26}$ because the area under $\textcolor{blue}{\sigma_L}\textcolor{red}{\pi_R}$ is 5 and the area under $\textcolor{red} {\pi_L}\textcolor{blue}{\sigma_R}$ is 26. 
 
\bigskip

\noindent {\bf References:}

\bigskip

\noindent [1] L. M. Butler, ``The $q$-log-concavity of $q$-binomial coefficients", {\em J. Combin. Theory Ser. A} {\bf 54} (1990), 54--63.

\medskip

\noindent [2] L. Carlitz and J. Riordan, ``Two element lattice permutations and their $q$-generalization", {\em Duke J. Math.} {\bf 31} (1964), 371--388.

\medskip

\noindent [3] J. F\"urlinger and J. Hofbauer, ``$q$-Catalan numbers",  {\em J. Combin. Theory Ser. A} {\bf 40} (1985), 248--264.

\medskip

\noindent [4] Li Liu and Yi Wang, ``On the log-convexity of combinatorial sequences'', to appear in {\em Advances in Applied Mathematics}, arXiv:math.CO/0602672.

\medskip

\noindent [5] P. A. MacMahon, {\em Combinatory Analysis}, Vol. I, Cambridge, 1915.
      \end{document}